\documentclass{article}

\usepackage{latexsym}

\usepackage{amsmath}

\usepackage{amsfonts}

\usepackage{amssymb}

\newcommand{\ms}{\medskip}

\newcommand{\bbc}{{\mathbb C}}

\newcommand{\bbr}{{\mathbb R}}

\newcommand{\co}{{\mathcal O}}

\newcommand{\fg}{{\mathfrak g}}

\newcommand{\dfg}{\fg^{\ast}}

\newcommand{\symm}{S(\fg)}

\newcommand{\symmn}{S_n(\fg)}

\newcommand{\poly}{P(\co)}

\newcommand{\cas}{C(\fg)}

\newcommand{\obid}{I(\co)}

\newcommand{\adl}{ad^{-1}(L)}

\newcommand{\jl}{J(L)}

\newcommand{\rdl}{r(J(L))}

\newcommand{\cpxob}{\co_{\bbc}}

\newcommand{\cpxalg}{\fg_{\bbc}}

\newcommand{\cpxgrp}{G_{\bbc}}

\newcommand{\cpxpoly}{\poly_{\bbc}}

\newcommand{\cpxid}{I_{\bbc}}

\newcommand{\realid}{K_{\bbr}}

\newcommand{\zeroset}{V(\cpxid)}

\newcommand{\pbrkt}{\{\; ,\,\}}

\newcommand{\pa}{Poisson algebra}

\newcommand{\proof}{\noindent{\it Proof.  }}

\newcommand{\de}{\textrm{d}}

\newcommand{\pl}{\partial}

\newcommand{\we}{\wedge}

\def\endproof{\hfill $\blacksquare$}

\newtheorem{thm}{Theorem}

\newtheorem{lem}{Lemma}

\newtheorem{prop}[thm]{Proposition}

\begin{document}



\title{POLYNOMIAL ALGEBRAS ON COADJOINT ORBITS OF SEMISIMPLE
LIE GROUPS}

\author{{\bf Mark J. Gotay}\thanks{Supported in part by NSF
grants 96-23083 and 00-72434. E-mail: gotay@math.hawaii.edu}
\\  Department of Mathematics \\ University of Hawai`i \\ 2565
The Mall \\ Honolulu, HI 96822  USA \\ 
\and {\bf Janusz Grabowski}\thanks{Supported by KBN, grant No.
2 P03A 031 17. E-mail: jagrab@mimuw.edu.pl} \\ Institute of
Mathematics
\\ University of Warsaw \\ ul. Banacha 2 \\ 02-097 Warsaw,
Poland \\
\and {\bf Bryon Kaneshige}\thanks{E-mail:
bryon@math.hawaii.edu} \\ Department of Mathematics \\ University of Hawai`i \\ 2565
The Mall \\ Honolulu, HI 96822  USA}

\date{August 19, 2000}

\maketitle





\begin{abstract}
We study the algebraic structure of the Poisson algebra $P(\co)$
of polynomials on a coadjoint orbit $\co$ of a 
semisimple Lie algebra. We prove that $\poly$ splits into a
direct sum of its center and its derived ideal. We also show
that $\poly$ is simple as a Poisson algebra iff $\co$ is
semisimple.
\end{abstract}





\begin{section}{Structure Theorems}

Let $\fg$ be a real (finite-dimensional) semisimple Lie algebra
with corresponding 1-connected Lie group $G$. 
It is well known that the dual space $\dfg$ carries the
structure of a linear Poisson manifold under the Lie-Poisson
bracket. The symplectic leaves of this Poisson structure are
the orbits of the coadjoint representation of $G$ on $\dfg$.

As the elements of $\fg$
may be regarded as linear functions on $\dfg$,
the symmetric algebra $\symm$ may be identified with the
algebra of polynomial functions on $\dfg$. Consequently, $\symm$
can be realized as a Poisson subalgebra of $C^\infty(\dfg)$.
(Equivalently, the Poisson bracket
$\pbrkt$ on $\symm$ can be obtained by setting $\{\xi,\eta\} =
[\xi,\eta]$ for $\xi, \eta \in \fg$ and extending to all of
$\symm$ via the Leibniz rule.)

Let $\symm'=\{\symm, \symm\}$ be the derived ideal, and let
$\cas$ denote the Lie center of the Poisson algebra
$\symm$. 

\begin{prop}

$\symm = \cas \oplus \symm'.$
\label{s}
\end{prop}

\proof  We have the decomposition 
$$\symm = \bigoplus^{\infty}_{n=0}\symmn$$ 
of the symmetric algebra into the finite-dimensional subspaces
$\symmn$  of homogeneous elements of
degree $n$. Each $\symmn$ is invariant with respect to the
adjoint action of $\fg$ on
$\symm$. Since every finite-dimensional representation of a
semisimple Lie algebra is completely reducible, it follows
that the
adjoint action of $\fg$ on
$\symm$ is itself completely reducible. According to
\cite[\S1.2.10]{Di}  we can then split
$$\symm = \cas \oplus \{\fg,\symm\}.$$ 
So we need only show that $\{\fg,\symm\}
= \symm'$.

Now, applying the identity
$$\{fg,h\} = \{f,gh\} + \{g,fh\}$$
to $f,g \in \fg$ and $h \in \symmn$, we see that
$\{S_2(\fg),S_n(\fg)\} \subset \{\fg,S_{n+1}(\fg)\}$.
Arguing recursively, we obtain 
$$\{S_m(\fg),S_n(\fg)\} \subset \{\fg,S_{n+m-1}(\fg)\},$$
from which the desired result follows.
\endproof

\ms

Let $\co$ be an orbit in $\dfg$. We can restrict polynomials on
$\dfg$ to functions on $\co$ thereby obtaining the (orbit)
polynomial algebra $\poly$ (which, however, may not be freely
generated as an associative algebra). We may identify $\poly$
with the quotient algebra $\symm / \obid$, where $\obid$ is the
ideal of elements vanishing on $\co$, with the canonical
projection
$$\rho_{\co} \colon \symm \to \symm / \obid \cong \poly.$$
Since $\co$ is a symplectic leaf of the Poisson structure on
$\dfg$, $\obid$ is a Lie ideal as well. 
Thus $\obid$ is a Poisson ideal (i.e., an associative
ideal which is also a Lie ideal) and hence
$\poly$ is a Poisson algebra of polynomial functions on the
symplectic leaf
$\co$. Note that since $\co$ is symplectic, the Lie center
$Z(P(\co)) = \bbr.$

We now show that the decomposition in Proposition~\ref{s} projects
to a similar decomposition of $\poly$.

\begin{thm}
$\poly = \bbr \oplus \poly'.$
\label{p}
\end{thm}

\proof  It is clear that $\cas$ projects onto constants on
$\co$ and $\rho_{\co}(\symm') = \poly'$, so that $\bbr +
\poly' = \poly$ by Proposition~\ref{s}. It remains to show that
$\poly' \cap \bbr = \{0\}$.

Now the restriction of the adjoint action of
$\fg$ on $\symm$ to the invariant subspace $\obid$ is
also completely reducible, so we can again use
\cite[\S1.2.10]{Di} to split
$\obid = \obid_{1} \oplus \obid_{2}$, where $\obid_{1} =
\obid \cap \cas$ and $\obid_{2} = \{\fg,\obid\} \subset \symm'$.

If $\poly' \cap \bbr \neq \{0\}$, then there is an
$f \in \obid$ such that $1+f \in \symm'$. Decomposing $f =
f_{1} + f_{2}$ with $f_{1} \in \obid_{1}$ and $f_{2} \in
\obid_{2}$, we get $(1+f_{1}) + f_{2} \in \symm'$, so $1+f_{1}
\in \cas \cap
\symm' = \{0\}$. Hence $f_{1} = -1$, and this contradicts the
fact that $f_{1} \in \obid_1 \subset \obid.$
\endproof

\ms Theorem~\ref{p} was already known in the case when $\co$ is
compact \cite{GGG}. In the $C^\infty$ context, one knows that if
$M$ is a compact symplectic manifold, then its \pa
$$C^\infty(M) = \bbr \oplus C^\infty(M)'$$
\cite{Av}, while if $M$ is noncompact
$$C^\infty(M) = C^\infty(M)'$$
\cite{Li}. Since in the smooth case $f\in C^\infty(M)'$ if 
and  only  if 
$f\eta$ is an exact  form,  where  $\eta$  is  the 
Liouville  volume  form, Theorem~\ref{p} thus suggests that the
polynomial Poisson (resp. de Rham) cohomology of a noncompact
coadjoint orbit
$\co$ may differ from its smooth Poisson (resp. de Rham)
cohomology. For example,  take 
$\co \subset sl(2,\bbr)^*$   to be the  one-sheeted 
hyperboloid 
$x^2+y^2-z^2=1$.  The  Poisson  tensor 
$$\Lambda=x\pl_y\we\pl_z+y\pl_z\we\pl_x - z\pl_x\we\pl_y$$ 
on $sl(2,\bbr)^*$ is polynomial and  the  induced 
symplectic  form 
$$\omega=x\,\de y\we\de z+y\,\de z\we\de x+z\,\de x\we\de y$$ 
on $\co$ is also  polynomial. As $\omega$ is a 
volume  form  on  the  non-compact manifold $\co$, it is
exact in the smooth category.  It is, however, not exact in the
polynomial category.  Indeed,  if  $\omega=\de\alpha$
for some polynomial 1-form  $\alpha$  on 
$\co$, then,  according  to  the  well-known  isomorphism 
between  Poisson and de Rham cohomology on a symplectic
manifold, we  would  have $[\Lambda, i_\alpha\Lambda]=\Lambda$,
where $[\ ,\ ]$ is  the  Schouten bracket and $[\Lambda,\cdot]$
is the  Poisson  cohomology  differential \cite{Va}. Writing
$\alpha=f\,\de x+g\,\de y+h\,\de z$,  where  $f,g,h$  are
polynomials,  this gives
$$\Lambda=H_x\we  H_f+H_y\we  H_g+H_z\we  H_h,$$ where $H_a$ is
the Hamiltonian vector field of  $a$.  Contracting $\Lambda$
with $\omega$ then yields
$1=\{x,f\}+\{y,g\}+\{z,h\}$---a contradiction with
Theorem~\ref{p}.
 
We remark that Theorem~\ref{p} need not hold if $\fg$ is not
semisimple. For instance, $\bbr^{2n}$ with its standard
symplectic structure is a coadjoint orbit of the Heisenberg
group H($2n$), but in this case $P(\bbr^{2n}) =
P(\bbr^{2n})'$.

\end{section}





\begin{section}{A Characterization of $\poly$}

We call a Lie algebra $L$ \emph{essentially simple} if
every Lie ideal of $L$ is either contained in the center $Z(L)$
of $L$ or contains the derived ideal $L'=[L,L]$. We say that a
\pa\ $P$ is \emph{simple} if the only Poisson ideals of $P$ are
$P$ and $\{0\}$.

\begin{prop}
Let $P$ be a unital \pa\ which has no nilpotent elements with
respect to the associative structure. If $P$ is simple, then it
is essentially  simple.
\label{es}
\end{prop}

\proof In view of \cite[Thm.1.10]{Gr}, if $L$ is a Lie ideal of
a unital \pa\ $P$ then 
\begin{equation}
\{P , \adl \} \subset \rdl,
\label{r}
\end{equation}
where
$$\adl = \{f \,|\, \{f, P \} \subset L \},$$
$\jl$ is the
largest associative ideal of $P$ contained in $\adl$, and
$\rdl$ is its radical, $$\rdl = \{f \,|\, f^n \in \jl \textrm{
for some } n=1,2,\ldots \}.$$
We recall from \cite[Thm. 1.6]{Gr} that $\jl$ is in fact a
Poisson ideal of $P.$

Suppose that $P' \not \subset L$.
Then $\adl \neq P$, so $\jl \neq P$, and thus
$\jl = \{0\}$  as $P$ is simple. Then $\rdl =
\{0\}$ since by assumption $P$ has no associative
nilpotents. Then (\ref{r}) gives
$$\{P, L\} \subset \{P, \adl \}  = \{0\},$$
i.e., $L \subset Z (P)$.
\endproof

\ms

In particular, the hypotheses of Proposition~\ref{es} are
satisfied by the polynomial algebra $\poly$. We now use this
Proposition to prove our main result.

\begin{thm}

The Lie algebra $\poly$ is essentially simple iff the orbit
$\co$ is semi\-simple.

\end{thm}

\proof  $(\Leftarrow)$ Assume $\co$ is semisimple and let
$\cpxob$ be the complexification of $\co$, i.e., the orbit in
$\cpxalg^{\:\:\ast}$ with respect to the complexified Lie group
$\cpxgrp$ which contains $\co$ in its real part. It
is well known that
$\cpxob$ is semisimple and that $\cpxob$
is an algebraic set in $\cpxalg^{\:\:\ast}$ \cite[\S3.8]{Ko}. If
$\poly$ were not essentially simple, then by
Proposition~\ref{es} we would have a proper Poisson ideal $I$ in
$\poly$, and so, after complexification, a proper  Poisson ideal
$\cpxid$ in $\cpxpoly := P_{\bbc}(\cpxob)$.

Let $\zeroset$ be the set of zeros of $\cpxid$ in $\cpxob$.
Since $\cpxob$ is algebraic, $\zeroset \neq \emptyset$, and
since $\cpxid$ is a Lie ideal, $\zeroset$ is
$\cpxgrp$-invariant and hence consists of orbits. This forces
$\zeroset = \cpxob$ and so $\cpxid = \{0\}$. Hence we have a
contradiction, since $\cpxid$ is proper.

$(\Rightarrow)$ Assume that $\co$ is not semisimple.
Complexifying as before, we get the complexified orbit
$\cpxob$ which is not semisimple. Now there exists a
semisimple orbit $\mathcal S$ in the Zariski closure of
$\cpxob$ \cite[\S3.8]{Ko}. Consider the Poisson ideal $K$ of
elements of
$\poly$ which vanish on $\mathcal S$. We claim that this ideal
is proper. Indeed, $K = \{0\}$ implies that $I(\mathcal S) =
I(\cpxob)$ ,whence $\mathcal S = {\rm cl}(\cpxob)$ as
$\mathcal S$  is an algebraic set. But this is impossible as
$\mathcal S$ and $\cpxob$ are distinct orbits. As well,
$K = \cpxpoly$ is impossible
as $\mathcal S \neq \emptyset.$

Now we will show that the existence of the proper Poisson ideal
$K$ in the complex Poisson algebra $\cpxpoly$ implies the
existence of a proper Poisson ideal $I$ in $\poly$. First, put
$$\realid = \{f \in \poly \,|\, f+ig \in K \textrm{ for some
}g \in \poly \}.$$ 
Since for $h \in \poly$, $f+ig \in K$ implies $(hf) +
i(hg) \in K$ and $\{h,f \} + i\{h,g \} \in K$, $\realid$ is
a Poisson ideal of $\poly$. Clearly 
$K \subset \realid + i\realid$ so that if $\realid = \{0\}$
then $K = \{0\}$. We can thus take $I = \realid$ as long as
$\realid \neq \poly$. 

If $\realid = \poly$, then there is $g
\in \poly$ such that $1 + ig \in K$.
Let 
$$K_0 = \{f \in \poly \,|\, if \in K \}.$$ 
Similarly as for
$\realid$, we can prove that $K_0$ is a Poisson ideal. Now $K_0
\neq \poly$, for otherwise
$K = \cpxpoly$.  We can then take $I = K_0$ provided $K_0
\neq \{0\}$. But in fact $K_0 \neq \{0\}$: Since
$$\{\poly, 1+ig\} = i\{\poly, g\}
\subset K,$$ 
$\{\poly, g\} \subset K_0$, and so $K_0 = \{0\}$
implies that $g \in Z(\poly) = \bbr$. So $1+ig \in K$ is a
nonzero constant, whence again $K = \cpxpoly$.

In any eventuality, we  now have a proper Poisson ideal $I$
of $\poly$. Of course, $I \not \subset Z(\poly) = \bbr$.
However, it may happen that $I \supset \poly'$, in which
case Theorem~\ref{p} forces $I = \poly'$. In this circumstance
we pass to the associative ideal $I^2$. Since
$$\{\poly, I^2\} \subset \{\poly , I\} I \subset I^2,$$ 
$I^2$ is also a Lie ideal. If $I^2 \neq I$, then $I^2$ is a
proper Lie ideal which  neither is contained  in  the center
nor contains the derived ideal.

To see that $I^2 \neq I$ for $I$ proper, we can use the
following.

\begin{lem}

If $P$ is a commutative unital ring with no zero divisors and
$I$ is a proper ideal which is finitely generated, then $I^2
\neq I$.

\end{lem}

\proof Assume that $x_1,\ldots,x_n$ are generators of $I$ and
$I^2 = I$. Then $x_i = \sum_{j=1}^n a_{ij}x_j$ for some $a_{ij}
\in I$, so that $\sum_{j=1}^n b_{ij}x_j = 0$ where $b_{ij} =
\delta_{ij} -a_{ij}$. Setting $B = \left(b_{ij}\right)$,
Cramer's Rule gives $x_i \det B = 0$ whence $\det B = 0$. But
$\det B \in \{1\} + I$ so $\det B \neq 0$.
\hfill $\blacktriangledown$

\ms

Thus $\poly$ is not essentially simple.
\endproof

\ms

The last part of this proof provides a converse to
Proposition~\ref{es} when $P = \poly$. In particular, we 
conclude that $\poly$ is simple if and only if $\co$ is
semisimple.

One can also see explicitly that $\poly$ is not essentially
simple when $\co$ is nilpotent as follows. 
Since a nilpotent orbit is conical \cite{Br}, 
$I(\co)$ is a homogeneous ideal. As a consequence, the notion of
homogeneous polynomial makes sense in
$P(\co)$. Let $P_k(\co)$ denote the subspace consisting
of all homogeneous polynomials of degree $k$. By virtue of
the commutation relations of $\fg$, 
$$\{P_k(\co),P_l(\co)\} \subset
P_{k+l-1}(\co),$$ 
whence each
$P_{(k)}(\co) = \oplus_{\ell \geq k}P_\ell(\co)$ for $k \geq 1$
is a proper Poisson ideal of $\poly$. 

\end{section}








\begin{thebibliography}{99}

\bibitem[Av]{Av} A. Avez, \emph{Remarques sur les
automorphismes infinit\'esimaux des vari\'et\'es symplectiques
compactes,} Rend. Sem. Mat. Univers. Politecn. Torino
{\bf 33} (1974--1975), 5--12.


\bibitem[Br]{Br} R. Brylinski,
\emph{Geometric quantization of nilpotent orbits,} J. Diff.
Geom. Appl. {\bf 9} (1998), 5--58.

\bibitem[Di]{Di} J. Dixmier, {\it
Enveloping Algebras,}  North-Holland, Amsterdam (1977). 

\bibitem[GGG]{GGG} M. J. Gotay, J. Grabowski, and H. B. 
Grundling,
\emph{An obstruction to quantizing compact symplectic
manifolds,} Proc. Amer. Math. Soc.
{\bf 28} (2000), 237--243.

\bibitem[Gr]{Gr} J. Grabowski, \emph{The Lie structure of
$\,C^*$ and Poisson algebras}, Studia Math. {\bf 81} (1985),
259--270.

\bibitem[Ko]{Ko} B. Kostant,
\emph{Lie group representations on polynomial rings,} Amer. J.
Math. {\bf 85} (1963), 327--404.

\bibitem[Li]{Li} A. Lichnerowicz,
\emph{Sur l'alg\`ebre des automorphismes infinit\'esimaux
d'une vari\'et\'e symplectique,} J. Diff. Geom. {\bf 9}
(1974), 1--40.

\bibitem[Va]{Va} I. Vaisman, {\it
Lectures on the Geometry of Poisson Manifolds,} Progress in
Math. {\bf 118}, Birkh\"auser, Boston (1994). 


\end{thebibliography}
\end{document}